\documentclass[10pt,a4paper]{article}
\usepackage[latin1]{inputenc}
\usepackage{amsmath}
\usepackage{amsfonts}
\usepackage{amssymb}
\usepackage{theorem}
\usepackage{mathrsfs}

\newtheorem{lemma}[subsection]{Lemma}
\newtheorem{corollary}[subsection]{Corollary}
\newtheorem{proposition}[subsection]{Proposition}
\newtheorem{example}[subsection]{Example}
\newtheorem{definition}[subsection]{Definition}
\newtheorem{remark}[subsection]{Remark}
\newtheorem{theorem}[subsection]{Theorem}

\title{Common Coupled Fixed Point Theorems for Maps in a $T_0$-ultra-quasi-metric Space}
\author{  Collins Amburo Agyingi,  \\ Ya\'e Ulrich Gaba,  \\ Department of Mathematics and Applied Mathematics,\\
University of Cape Town, Rondebosch 7701, South Africa}
\begin{document}

\maketitle

\begin{abstract}
In this article, we prove the existence of common fixed points for a pair of maps on a $q$-spherically complete $T_0$-ultra-quasi-metric space. The present article is a generalization, in the assymmetric setting of the paper of Rao et al.\cite{Rao}. The key point in the proof is the use of Zorn's lemma. We construct an appropriate chain and show that it has a maximal element, from which we extract the fixed point we are looking for. The choice of the sets, here open balls, is characteristic of this type of problems and the contraction condition are essential, specially when we are to establish the uniqueness of the fixed point. \\
Keywords-\textit{$q$-spherical completeness; $T_0$-ultra-quasi-metric space; common fixed point.}
\end{abstract}

\section{\textsc{INTRODUCTION}}
To prove results about fixed points or common fixed points for maps in $T_0$-ultra-quasi-metric space satisfying some strict contractive conditions, supcontinuity of the maps is assumed as well as the space is required to be joincompact. It turns out that this is not the case for $q$-spherically complete $T_0$-ultra-quasi-metric space. In this situation, supcontinuity of maps is not necessary to obtain fixed points. A fixed points theorem for maps in $T_0$-ultra-quasi-metric space had been proved recently (compare \cite{Agyingifixed}) and the result has been extended in \cite{Agyingiproximity} for multi-valued maps.

For recent results in the area of Asymmetric Topology, the reader is adviced to consult \cite{Agyingitight, Agyingiendpoint2}.

\section{\textsc{PRELIMINARIES}}
\definition (Compare \cite[page 2]{KunziOlela}) Let $X$ be a set and $d:X\times X \rightarrow [0,\infty)$ be a function mapping into the set $[0,\infty)$ of non-negative reals. Then $d$ is an {\em ultra-quasi-pseudometric} on $X$ if 

(a) $d(x,x)=0$ for all $x\in X$, and 

(b) $d(x,z)\leq max\lbrace d(x,y),d(y,z)\rbrace$ whenever $x,y,z\in X$.

The conjugate $d^{-1}$ of $d$ where $d^{-1}(x,y)=d(y,x)$ whenever $x,y\in X$ is also an ultra-quasi-pseudometric on $X$.  \\
If $d$ also satisfies the following condition (known as the $T_0$-condition):

(c) for any $x,y\in X, \; d(x,y)=0=d(y,x)$ implies that $x=y$, then $d$ is called a {\em $T_0$-ultra-quasi-metric} on $X$. Notice that $d^s=\sup\lbrace d,d^{-1}\rbrace= d\vee d^{-1}$ is an {\em ultra metric} on $X$ whenever $d$ is a $T_0$-ultra-quasi-metric on $X$  .

In the literature, $T_0$-ultra-quasi-metric spaces are also know as non Archimedean $T_0$-quasi-metric spaces. The set of open balls $\lbrace \lbrace y\in X:d(x,y)<\epsilon\rbrace:x\in X, \epsilon>0\rbrace$ yields a base for the topology $\tau(d)$ induced by $d$ on $X$.

\example (Compare \cite[Example 3]{Bonsangue}) \label{Bon} Let $X=[0,\infty)$. Define for each $x,y\in X$, $n(x,y)=x$ if $x>y$, and $n(x,y)=0$ if $x\leq y$. It is not difficult to check that $(X,n)$ is a $T_0$-ultra-quasi-metric space.

Notice also that for $x,y\in [0,\infty)$, we have $n^s(x,y)=\max\lbrace x,y \rbrace$ if $x\neq y$ and $n^s(x,y)=0$ if $x=y$. The ultra metric $n^s$ is complete on $[0,\infty)$ since $n$ and $n^{-1}$ are complete on $[0,\infty)$ (compare \cite[Example 2]{KunziOlela}).

\definition \label{supcontinuous} {Let $(X,d)$ and $(Y,m)$ be two ultra-quasi-pseudometric spaces. We say that a map $f:(X,d)\to (Y,m)$ is {\em supcontinuous} if $f:(X,d^s)\to (Y,m^s)$ is continuous.

Let $(X,d)$ be an ultra-quasi-pseudometric space. Let $x\in X$ and \\$r\in [0,\infty)$. By $C_d(x,r)$ we mean the closed ball \[ C_d(x,r)=\lbrace y\in X: d(x,y)\leq r\rbrace \] of radius $r$ centred at $x$.

\lemma (Compare \cite[Lemma 9]{KunziOlela}) If $(X,d)$ is an ultra-quasi-pseudometric space and $x,y\in X$ and $r,s\in [0,\infty)$, then we have that \[ C_d(x,r)\cap C_ {d^{-1}}(y,s)\neq \emptyset \] if and only if \[d(x,y)\leq \max\{ r,s\}.\]

\definition (Compare \cite[Definition 2]{KunziOlela}) Let $(X,d)$ be an ultra-quasi-pseudometric space. Let $(x_i)_{i\in I}$ be a family of points in $X$ and let $(r_i)_{i\in I}$ and $(s_i)_{i\in I}$ be families of non-negative real numbers. We shall say that the family $(C_d(x_i,r_i),C_{d^{-1}}(x_i,s_i))_{i\in I}$ has the mixed binary intersection property provided that \[ d(x_i,x_j)\leq \max\lbrace r_i, s_j\rbrace \]whenever $i,j\in I$.

\definition We say that $(X,d)$ is $q$-spherically complete provided that each family $(C_d(x_i,r_i),C_{d^{-1}}(x_i,s_i))_{i\in I}$ possessing the mixed binary intersection property also satisfies  $$ \bigcap_{i\in I}\left( C_d(x_i,r_i)\cap C_{d^{-1}}(x_i,s_i) \right)\neq \emptyset.$$

For an example of a $q$-spherically complete ultra-quasi-metric space, the reader is adviced to check \cite[Example 2]{KunziOlela}.

\proposition (Compare \cite[Proposition 2]{KunziOlela})

(a) Let $(X,d)$ be an ultra-quasi-pseudometric space. Then $(X,d)$ is $q$-spherically complete if and only if $(X,d^{-1})$ is $q$-spherically complete.

(b) Let $(X,d)$ be a $T_0$-ultra-quasi-metric space. If $(X,d)$ is $q$-spherically complete, then $(X,d^s)$ is spherically complete.

\definition An ultra-quasi-pseudometric space $(X,d)$ is called bicomplete provided that the ultra-pseudometric $d^s$ on $X$ is complete.

\proposition (Compare \cite[Proposition 3]{KunziOlela}) Each $q$-spherically complete $T_0$-ultra-quasi-metric space $(X,d)$ is bicomplete.

In \cite{Agyingifixed}, Agyingi proved the following:

\theorem Let $(X,d)$ be a $q$-spherically complete $T_0$-ultra-quasi-metric space. If $T:X\to X$ is a mapping such that 
$$ d(Tx,Ty)< \max\{ d(x,y),d(Tx,x),d(y,Ty)\} ,$$
for all $x,y\in X$, with $0<\min\{d(x,y),d(y,x)\}$, then $T$ has a unique fixed point.

Now, we extend this theorem for a pair of maps of Jungck type.
\section{\textsc{Main Results}}
\definition Let $(X,d)$ be an ultra-quasi-pseudometric space, $f:X\to X$ and $T:X\to X$ be two self maps on $X$. We say that $f$ and $T$ {\em coincidentally commute} at $z\in X$ if $Tfz=fTz$.

\theorem\label{mainresult}  Let $(X,d)$ be a $q$-spherically complete $T_0$-ultra-quasi-metric space. If $f$ and $T$ are self map on $X$ satisfying 
\begin{equation}
T(X)\subseteq f(X),
\end{equation}
\begin{equation}\label{2}
d(Tx,Ty)< \max\{ d(fx,fy),d(Tx,fx),d(fy,Ty)\} ,
\end{equation}
for all $x,y\in X$, with $0<\min\{d(x,y),d(y,x)\}$, then there exists $z\in X$ such that $fz=Tz$. 

Further, if $f$ and $T$ coincidentally commute at $z$, then $z$ is the unique fixed point of $f$ and $T$.

{\bf Proof.}

Let $a\in X$. Let us denote by 

$$ C_d^a = C_{d}(fa,\alpha_a) \text{ and }    C_{d^{-1}}^a= C_{d^{-1}}(fa,\alpha_a)                                   ,$$
with $\alpha_a:= d^s(fa,Ta)=\max\{d(Ta,fa),d(fa,Ta)\}$ and set 
$$C^a = C_d^a \cap  C_{d^{-1}}^a .$$

Let $\mathcal{A}:=\lbrace C^a: a\in X \rbrace$. Define the relation $\preceq$ on $\mathcal{A}$ by 
   $$ C^a \preceq C^b \text{ if and only if }  C^b \subseteq C^a .$$
   
  Then $(\mathcal{A},\preceq)$ is a partially ordered set. Since the verification of this fact is trivial, we leave it to the interested reader.

Let $\mathcal{A}_1$ be a nonempty chain in $\mathcal{A}$. Then by $q$-spherical completeness of $(X,d)$, we have that \[ \bigcap_{C^a\in \mathcal{A}_1}C^a=C\neq \emptyset.\]
Clearly, $C\cap f(X) \neq \emptyset$. It follows from the mature of the balls considered, the centers are of the form $fa, a\in X$. Moreover, we consider a totally ordered family of such balls. In fact, 
 \[ \underset{C^a \in \mathscr{A}_1}{\bigcap}C^a=C^y \] for some $y\in X$. For more clarity in this part, we refer the reader to \cite{Rao}.
 
Let $fb\in C$ and $C^a\in \mathcal{A}_1$. Then we have \[d(fa,fb)\leq \alpha_a \; \text{and}\; d(fb,fa)\leq \alpha_a. \]

Let now $x\in C^b$. Then \[d(fb,x)\leq \alpha_b \; \text{and}\; d(x,fb)\leq \alpha_b . \]

\begin{align*}
d(fb,x)&\leq \max\{d(Tb,fb),d(fb,Tb)\} \\
&\leq \max \lbrace d(Tb,Ta),d(Ta,fa),d(fa,fb), d(fb,fa), d(fa,Ta), d(Ta,Tb)\rbrace \\
&\leq \max\lbrace d(Tb,Ta),d(Ta,Tb),\alpha_a \rbrace \\
&\leq \max\lbrace d(Tb,fb),d(fb,fa),d(fa,Ta), d(Ta,fa),d(fa,fb),d(fb,Tb)\}\\
&\leq \alpha_a
\end{align*}

From the above inequality, we have now that

\begin{align*}
d(fa,x)&\leq \max\lbrace d(fa,fb),d(fb,x)\rbrace \\
&= \alpha_a
\end{align*}
which means that $x\in C_d(fa,\alpha_)$. We have thus shown that 
\begin{equation} \label{inclusion1} 
C_d(fb,\alpha_b)\subseteq C_d(fa,\alpha_a).
 \end{equation} 

By a similar computation, one can show that
\begin{equation} \label{inclusion2} 
C_{d^{-1}}(fb,\alpha_b)\subseteq C_{d^{-1}}(fa,\alpha_a). 
\end{equation}

By Equations \eqref{inclusion1} and \eqref{inclusion2}, we have that for all $C^a\in \mathcal{A}_1,\;C^b\subseteq C^a$, which means that $C^a\preceq C^b$ for all $C^a\in \mathcal{A}_1$. Thus $C^b$ is an upper bound in $\mathcal{A}$ for the chain $\mathcal{A}_1$. By Zorn's lemma we conclude that $\mathcal{A}$ has a maximal element, say, $C^u,\; u\in X$. We claim that $Tu=fu$.

Suppose on the contrary that $fu\neq Tu$. Since $Tu \in T(X)\subseteq f(X)$, there exists $w\in X$ such that $Tu=fw$. From \eqref{2}, we have
\begin{align*}
  d(fw,Tw)&=d(Tu,Tw) \\
           & < \max\{ d(fu,fw),d(Tu,fu),d(fw,Tw)  \} \\
           & < \max \{d(fu,fw),d(fw,fu),d(fw,Tw) \} \\
           & = \max \{d(fu,fw),d(fw,fu)\}.
\end{align*}

and 
\begin{align*}
  d(Tw,fw)&=d(Tw,Tu) \\
           & < \max\{ d(fw,fu),d(fu,Tu),d(Tw,fw)  \} \\
           & < \max \{d(fw,fu),d(fu,fw),d(fw,Tw) \} \\
           & = \max \{d(fu,fw),d(fw,fu)\}.
\end{align*}
Thus $fu \notin C^w$ and hence $C^u\subsetneq C^w$. It is a contradiction to the maximality of $C^u$. We therefore conclude that $Tu=fu.$
 
Moreover, if $f$ and $T$ coincidentally commute at $u$, then $f^2u=f(fu)=f(Tu)=T(fu)=T^2u.$
 
Suppose that $fu\neq u$. From condition \eqref{2}), we have that
 
 \begin{align*}
 d(Tfu,Tu) & < \max\{ d(f^2u ,fu), d(Tfu,f^2u),d(fu,Tu)\} \\
           & =  d(f^2u ,fu) = d(Tfu,Tu),
 \end{align*}
 since $d(Tfu,f^2u)=d(fu,Tu)=0$. 
 
Similarly we can prove that $d(Tu,Tfu)<d(Tu,Tfu)$. The above inequalities gives a contradiction and so we must have that $u=fu=Tu.$
 
Let us now prove uniqueness.\\
Suppose that there is another common fixed point, i.e., there exists $z\in X$ such that $f(z)=z=Tz$ and $z\neq u$. We shall examine two cases.\\

Case 1: Suppose $d(z,u)>0$. Then we have that 
\[d(z,u)=d(Tz,Tu)<\max\lbrace d(fz,fu),d(Tz,fz),d(fu,Tu)\rbrace=d(z,u),\]
which is a contradiction.\\

Case 2: Suppose now that $d(u,z)>0$. Then we get \[d(u,z)=d(Tu,Tz)<\max\lbrace d(fu,fz),d(u,z),d(fz,Tz)\rbrace=d(u,z),\]which is a contradiction. 

Thus we must have that $z=u$.

$\hfill{\Box}$

\corollary Theorem \ref{mainresult} holds if inequality (\ref{2}) is replaced by
\begin{equation}\label{3}
d(Tx,Ty)< \max\{ d(Tx,Ty),d(Tx,fx),d(fy,Ty),d(fx,Ty),d(Tx,fy)\},
\end{equation}
for all $x,y\in X$, with $0<\min\{d(x,y),d(y,x)\}$.

{\bf Proof.}

Since $d(fx,Ty) \leq \max\{  d(fx,fy),d(fy,Ty)\}$ and $d(Tx,fy) \leq \max\{  d(Tx,fx),d(fx,fy)\}$, hence inequality \eqref{3} implies inequality \eqref{2}.

$\hfill{\Box}$

\corollary Taking $f=I$ (the identity map), we obtain Theorem 1 of \cite{Agyingifixed}.

\section{\textsc{MULTI-VALUED MAPS}}

Now we generalize Theorem \ref{mainresult} when $T$ is a multi-valued map. 

Let $\mathcal{P}_0(X):=2^X \setminus \emptyset$ where $2^X$ denotes the power set of $X$. For $x\in X$ and $A,B \in \mathcal{P}_0(X)$, we set:

$$ d(x,A)= \inf\{ d(x,a),a\in A\}  \text{ and }  d(A,x)= \inf\{ d(a,x),a\in A\},$$
and define the Hausdorff ultra-quasi-pseudometric  $H$ by 

$$H(A,B)= \max \left\lbrace \underset{a\in A}{\sup}\ d(a,B), \underset{b\in B}{\sup} \ d(A,b)   \right\rbrace.$$

Then $H$ is an extended ultra-quasi-pseudometric on $\mathcal{P}_0(X)$. Moreover, we know from \cite{Kunzi} that the restriction of $H$ to $S_{cl}= \{A\subseteq X: A= cl_{\tau (d)}A \cap cl_{\tau (d^{-1})}A  \}$ is an extended $T_0$-ultra-quasi- metric.

\definition
For a non Archimedian $T_0$-quasi-metric space$(X,d)$, we denote by $2^X_j$ the space of all nonempty joincompact subsets in $X$ with the Hausdorff ultra-quasi-metric $H$.

We have the following result due to Agyingi \cite{Agyingifixed}.

\theorem (compare \cite{Agyingifixed}) Let $(X,d)$ be a $q$-spherically complete $T_0$-ultra-quasi-metric space. If $T:X\to 2^X_j$ is a mapping such that 
\[ H(Tx,Ty)< \max\{ d(x,y),d(Tx,x),d(y,Ty)\} \text{ for all }x,y\in X,\ x\neq y,\]
then $T$ has a unique fixed point, i.e there exists $x\in X$ such that $x\in Tx$.

\definition Let $(X,d)$ be an ultra-quasi-pseudometric space, $f:X\to X$ and $T:X\to  2^X_j $. $f$ and $T$ are said to be coincidentally commuting at $z\in X$ if $fz \in Tz $ implies that $fTz\subseteq Tfz$.

\theorem\label{multivalued} (compare \cite[Theorem 9]{Rao}) Let $(X,d)$ be a $q$-spherically complete $T_0$-ultra-quasi-metric space. Let $f:X\to X$ and $T:X\to  2^X_j $ be two maps satisfying 
\begin{equation}\label{6}
Tx\subseteq f(X) \quad \forall x\in X,
\end{equation}
\begin{equation}\label{7}
H(Tx,Ty)< \max\{ d(fx,fy),d(Tx,fx),d(fy,Ty)\},
\end{equation}
for all $x,y\in X$, with $0<\min\{d(x,y),d(y,x)\}$. Then there exists $z\in X$ such that $fz \in Tz$. 

Further, assume that 
\begin{equation}\label{8}
d^s(fu,fx) \leq \min\{H(Tu,Tfy), H(Tfy,Tu)\} \quad \forall x,y,u\in X \text{ with } fx\in Ty,
\end{equation}
and
 \begin{equation}\label{10}
 f \text{ and } T \text{ are coincidentally commuting at } z.
 \end{equation} 
Then $fz$ is the unique fixed point of $f$ and $T$.

{\bf Proof.}

Let $a\in X$ and denote by \[ C^d_a=C_d(fa,\beta_a)\; \text{ and }\ ; C^{d^{-1}}_a = C_{d^{-1}}(fa,d(fa,Ta)) \] the closed balls with centers at $fa\in X$ and 

radius $$\beta_a:=\max\{d(Ta,fa),d(fa,Ta)\}$$   with 

$$d(Ta,fa)=\inf\lbrace d(z,fa):z\in Ta\rbrace$$ and $$d(fa,Ta)=\inf\lbrace d(fa,z):z\in Ta\rbrace .$$

 Put\[ C_a= C^d_a\cap C^{d^{-1}}_a.\] Let $\mathcal{A}$ be the collection of all such closed balls $C_a$ such that $a$ runs over $X$. Define $\preceq$ on $\mathcal{A}$ by \[ C_a\preceq C_b\; \text{if and only if}\; C_b\subseteq C_a.\] Then $(\mathcal{A},\preceq)$ is a partially ordered set. We leave the verification of this fact to the interested reader.

Let $\mathcal{A}_1$ be a nonempty chain in $\mathcal{A}$. Then by $q$-spherical completeness of $(X,d)$, we have that \[ \bigcap_{C_a\in \mathcal{A}_1}C_a=C\neq \emptyset.\]
Here again , it is clear that  $C\cap f(X) \neq \emptyset$.

Let $fb\in C$ and $C_a\in \mathcal{A}_1$. Then we have \[d(fa,fb)\leq \beta_a \; \text{and}\; d(fb,fa)\leq \beta_a. \]

Let us choose $u\in Ta$ such that $d^s(fa,u)=d^s(fa,Ta).$ Notice that this is possible since the map $:Ta\to \mathbb{R}$ defined by $u\mapsto d(x,u)$ is uniformly continuous with respect to the usual metric on $\mathbb{R}$

With $u\in Ta$ satisfying the above condition and $fb\in C$, we have

\begin{align*}
d(fb,Tb)&= \inf\lbrace d(fb,c):c\in Tb\rbrace \\
&\leq \max\lbrace d(fb,fa),d(fa,u),\inf\lbrace d(u,c):c\in Tb\rbrace \rbrace \\
&\leq \max\lbrace \beta_a, H(Ta,Tb) \rbrace \\
&< \max\lbrace \beta_a, d(Ta,fa),d(fb,Tb),d(fa,fb) \rbrace \\
&= \max\lbrace \beta_a, d(fb,Tb)\rbrace \\
\end{align*}
which is possible only when $d(fb,Tb)<\beta_a.$

By a similar computation, we have that $d(Tb,fb)<\max\lbrace  \beta_a,d(Tb,fb)\rbrace$ which is possible only when $d(Tb,fb)< \beta_a.$

Let now $x\in C_b$. Then \[d(fb,x)\leq  \beta_b< \beta_a \ ; \text{and}\; d(x,fb)\leq \beta_b< \beta_a.\]

We have now that:
\begin{align*}
d(fa,x)&\leq \max\lbrace d(fa,fb),d(fb,x)\rbrace \\
&\leq \beta_a
\end{align*}
which means that $x\in C_d(fa,\beta_a)$. We have thus shown that
 \begin{equation} \label{inclusion11} 
 C_d(fb,\beta_b)\subseteq C_d(fa,\beta_a).
 \end{equation}

Similarly we can show that 
\begin{equation} \label{inclusion22} 
C_{d^{-1}}(fb,\beta_b) \subseteq C_{d^{-1}}(fa,\beta_a). 
\end{equation}

Equations \eqref{inclusion11} and \eqref{inclusion22} imply that for all $C_a\in \mathcal{A}_1,\;C_b\subseteq C_a$. In other words, this says that $C_a\preceq C_b$ for all $C_a\in \mathcal{A}_1$. Thus $C_b$ is an upper bound in $\mathcal{A}$ for the chain $\mathcal{A}_1$. We therefore conclude by Zorn's lemma that $\mathcal{A}$ has a maximal element, say, $C_z,\; z\in X$. We shall prove that $fz\in Tz$. We do this by contradiction.

Suppose on the contrary that $fz\not \in Tz$. Then from \eqref{6} there exists $fz^{\ast}\in Tz,\; fz^{\ast}\neq fz,$ such that $d^s(fz,fz^{\ast})=d^s(fz,Tz).$
\begin{align*}
d(fz^{\ast},Tz^{\ast})&\leq H(Tz,Tz^{\ast})\\
&< \max\lbrace d(fz,fz^{\ast}),d(Tz,fz),d(fz^{\ast},Tz^{\ast})\rbrace \\
& \leq \{\beta_z,d(fz^{\ast},Tz^{\ast})\}
\end{align*}
which is possible only if $d(fz^{\ast},Tz^{\ast})<\beta_z.$ 

Similarly we have that $d(Tz^{\ast},fz^{\ast})<\max\lbrace \beta_z, d(Tz^{\ast},fz^{\ast})$ which is possible only if $d(Tz^{\ast},fz^{\ast})<\beta_z$.

Let $y\in C_{z^{\ast}}$. 
Then \[ d(y,fz^{\ast})\leq \beta_{z^*}< \beta_z \ ; \text{and}\; d(fz^{\ast},y)\leq \beta_{z^*}< \beta_z. \]

Also \[ d(y,fz)\leq \max\lbrace d(y,fz^{\ast}),d(fz^{\ast},fz)\rbrace \leq \beta_z \] and \[ d(fz,y)\leq \max \lbrace d(fz,fz^{\ast}),d(fz^{\ast},y)\rbrace \leq \beta_z.\] This just means that $y$ belongs to the ball $C_{z}$, so that $C_{z^{\ast}}\subseteq C_{z}.$ 

The inequality $d^s(fz,fz^{\ast})=d^s(fz,Tz)>d^s(fz^{\ast},Tz^{\ast})$ implies that $z$ does not belong to the ball $C_{z^{\ast}}$ and this implies that $C_{z^{\ast}} \subsetneq C_z$ which contradicts the maximality of $C_z$. So we must have that $fz\in Tz$.

Further, assume \eqref{8}, and \eqref{10}. Form \eqref{8}, we have

\[
d(fz,f^2z)\leq H(Tfz,Tfz)= 0, 
\]
and 
\[
d(f^2z,f^2z) \leq H(Tfz,Tfz)= 0.
\]
This means that $d(fz,f^2z)=0=d(f^2z,f^2z)$, hence $f^2z=fz$.
From \eqref{10}, $fz=f^2z \in fTz \subseteq Tfz$. Thus $fz$ is a common fixed point of $f$ and $T$.

We shall now prove uniqueness.\\
Suppose there is $z^*$ such that $fz^*\neq fz$ and $fz^*=f^2z^* \in Tfz^*$. From \eqref{7}, \eqref{8}, we have

\begin{align*}
 d(fz,fz^*)=d(f^2z,f^2z^*) &\leq H(Tf^2z,Tfz^*)\\
            &\leq H(Tfz,Tfz^*) \\
            & < \max \{ d(f^2z,f^2z^*),d(f^2z,Tfz),d(f^2z^*,Tfz^*)\} \\
            & = d(fz,fz^*)
\end{align*}

and 

\begin{align*}
 d(fz^*,fz)&=d(f^2z^*,f^2z) \leq H(Tf^2z^*,Tfz)\\
            &\leq H(Tfz^*,Tfz) \\
            & < \max \{ d(f^2z^*,f^2z),d(f^2z^*,Tfz^*),d(f^2z,Tfz)\} \\
            & = d(fz^*,fz).
\end{align*}
This implies that $fz^*= fz$. Thus $z=f^2z$ is the unique common fixed point of $f$ and $T$.

$\hfill{\Box}$

\remark If $f=I$ (Identity map), then the first part of Theorem \ref{multivalued} is the main theorem of Agyingi \cite{Agyingiproximity}.}

\end{document}